\documentclass[11pt]{amsart}
\usepackage{amsmath}
\usepackage{amsfonts}
\usepackage{amsthm}
\usepackage{pictex}
\usepackage{graphicx}

\def\diag{\mathop{\mathrm{diag}}\nolimits}
\def\span{\mathop{\mathrm{span}}\nolimits}

\newcommand{\Tr}{{{\operatorname{Tr}}}}
\newtheorem{theorem}{Theorem}[section]

\newtheorem{Remark}[theorem]{Remark}
\newenvironment{remark}{\begin{Remark}\rm}{\end{Remark}}
\numberwithin{equation}{section}

\begin{document}
\sc
\title[Multiple Hermite and Laguerre polynomials]
{Integral representations for multiple Hermite and multiple
Laguerre polynomials}

\author{Pavel M. Bleher}
\address{Department of Mathematical Sciences, Indiana
University-Purdue University Indianapolis,
402 N. Blackford St., Indianapolis, IN 46202, U.S.A.}
\email{bleher@math.iupui.edu}
\author{Arno B.J. Kuijlaars}
\address{Department of Mathematics,
Katholieke Universiteit Leuven,
Celestijnenlaan 200 B,
B-3001 Leuven
BELGIUM}
\email{arno@wis.kuleuven.ac.be}

\dedicatory{Dedicated to Pierre van Moerbeke on his sixtieth
birthday}
\date{\today}

\thanks{The first author was supported in part by
NSF Grant DMS-9970625. The second author was supported by projects
G.0176.02 and G.0455.04 of FWO-Flanders.}

\rm

\begin{abstract}
We give integral representations for multiple Hermite and
multiple Hermite polynomials of both type I and II. We also
show how these are connected with double integral representations
of certain kernels from random matrix theory.
\end{abstract}

\maketitle

\section{Multiple orthogonal polynomials}

Multiple orthogonal polynomials are an extension of  orthogonal
polynomials that play a role in the random matrix ensemble with an
external source
\begin{equation} \label{ensemble}
    \frac{1}{Z_n} e^{- \Tr (V(M) - AM)} dM
\end{equation}
defined on $n \times n$ Hermitian matrices, see \cite{BK1,BK2,DK}.
Here $A$ is a fixed $n \times n$ Hermitian  matrix and $V : \mathbb R \to \mathbb R$ is
a function with enough increase at $\pm \infty$ so that the
integral
\[ Z_n = \int e^{- \Tr (V(M) - AM)} dM \]
converges. Random matrices with external source were introduced and studied
by Br\'ezin and Hikami \cite{BH1,BH2,BH3,BH4,BH5}, and
P. Zinn-Justin \cite{ZJ1,ZJ2}.

In what follows, we assume that $A$ has $m$ distinct eigenvalues
$a_1, \ldots, a_m$ of multiplicities $n_1, \ldots, n_m$. We consider
$m$ fixed and use  multi-index notation $\vec{n} = (n_1,
\ldots, n_m)$ and $|\vec{n}| = n_1 + \cdots + n_m$.

The average characteristic polynomial $P_{\vec{n}}(x) = \mathbb E \left[\det(xI - M) \right]$
of the ensemble (\ref{ensemble}) is a monic polynomial of degree $|\vec{n}|$
which  satisfies for $k=1, \ldots, m$,
\begin{equation} \label{typeII}
    \int_{-\infty}^{\infty} P_{\vec{n}}(x) x^j w_k(x) dx  = 0,
        \qquad j=0, \ldots, n_k-1,
\end{equation}
where
\begin{equation} \label{weights}
    w_k(x) = e^{- (V(x) - a_k x)},
\end{equation}
see \cite{BK1}.
The relations (\ref{typeII}) characterize the polynomial $P_{\vec{n}}$ uniquely.
For $A = 0$, we have that $P_{\vec{n}}$ is the usual orthogonal polynomial with respect
to the weight $e^{-(V(x))}$, which is a well-known fact from random
matrix theory. For general  $m$, the relations (\ref{typeII}) are multiple
orthogonality relations with respect to the weights (\ref{weights})
and the polynomial $P_{\vec{n}}$ is called a multiple orthogonal polynomial of type II.

The multiple orthogonal polynomials of type I consist of a vector
\begin{equation} \label{degreeAk}
    (A_{\vec{n}}^{(1)}, A_{\vec{n}}^{(2)}, \ldots, A_{\vec{n}}^{(m)}),
    \qquad \deg A_{\vec{n}}^{(k)} \leq n_k - 1,
\end{equation}
of polynomials such that the function
\begin{equation} \label{Qfunction}
    Q_{\vec{n}}(x) = \sum_{k=1}^m A_{\vec{n}}^{(k)}(x) w_k(x)
\end{equation}
satisfies
\begin{equation} \label{typeI}
    \int_{-\infty}^{\infty} x^j Q_{\vec{n}}(x) dx = \left\{ \begin{array}{ll}
        0, & \quad j = 0, \ldots, |\vec{n}|-2, \\
        1, & \quad j = |\vec{n}|-1.
        \end{array} \right.
\end{equation}
The polynomials $A_{\vec{n}}^{(k)}$ are uniquely determined by the
degree requirements (\ref{degreeAk}) and the type I
orthogonal relations (\ref{typeI}).

By the Weyl integration formula, the random matrix ensemble (\ref{ensemble})
has the following joint eigenvalue distribution
\[ \frac{1}{\tilde{Z}_n}
    \prod_{j=1}^n e^{- V(\lambda_j)} \left( \int e^{A U \Lambda U^*} dU \right)
        \prod_{j <k} (\lambda_j - \lambda_k)^2 \ d\lambda_1 d\lambda_2 \cdots d\lambda_n \]
where $dU$ is the normalized Haar measure on the unitary group $U(n)$
and $\Lambda = \diag(\lambda_1, \ldots \lambda_n)$. Using
the confluent form of the Harish-Chandra / Itzykson-Zuber formula
\cite{HC,IZ} to evaluate
the integral $\int e^{A U \Lambda U^*} dU$, we find that the
joint eigenvalue distribution can be written as
\begin{equation} \label{biorthogonal1}
    \frac{1}{Z'_n} \det\left(f_j(\lambda_k)\right)_{1\leq j, k \leq n}
    \det\left(g_j(\lambda_k)\right)_{1\leq j,k \leq n} \ d\lambda_1 d\lambda_2 \cdots d\lambda_n,
\end{equation}
where $f_j(x) = x^{j-1}$ and $g_1, \ldots, g_n$ are the functions
\begin{equation} \label{basisfunctions}
    x^j e^{-(V(x)-a_kx)}, \qquad k = 1, \ldots, m, \quad j=0,1, \ldots, n_k-1
\end{equation}
taken in some arbitrary order. Then (\ref{biorthogonal1}) is a
biorthogonal ensemble in the sense of Borodin \cite{Bor}, which in
particular implies that the eigenvalue point process is determinantal,
that is, there is a kernel $K(x,y)$ such that
the $k$ point correlation functions have determinantal form
\begin{equation} \label{correlations}
    \det\left( K(\lambda_i, \lambda_j) \right)_{1\leq i,j\leq k}.
\end{equation}

The multiple orthogonal polynomials of types II and I can be used to
biorthogonalize the functions $f_j$ and $g_j$ and to give an explicit
formula for $K$. Choose a sequence of
multi-indices $\vec{n}_0, \vec{n}_1, \ldots, \vec{n}_n = \vec{n}$
such that $|\vec{n}_j| = j$ and $\vec{n}_{j} \leq \vec{n}_{j+1}$
and define
\[ p_j(x) = P_{\vec{n}_j}(x), \qquad q_j(y) = Q_{\vec{n}_{j+1}}. \]
Then $p_{j+1} \in \span\{f_1, \ldots, f_j\}$, and
$q_{j+1} \in \span\{g_1, \ldots, g_j\}$ for a suitable ordering
of the functions (\ref{basisfunctions}).
In addition we have the biorthogonality
\[ \int_{-\infty}^{\infty} p_i(x) q_j(x) dx = \delta_{i,j}. \]
As in \cite{Bor} it then follows that
\begin{equation} \label{kernelK}
    K(x,y) = \sum_{j=0}^{n-1} p_j(x) q_j(y).
\end{equation}
By the Christoffel-Darboux formula for multiple orthogonal
polynomials \cite{DK} the kernel $K$ satisfies
\begin{equation} \label{CDkernel}
    (x-y) K(x,y) =
    P_{\vec{n}}(x) Q_{\vec{n}}(y)
    - \sum_{k=1}^m \frac{h_{\vec{n}}^{(k)}}{h_{\vec{n}-\vec{e}_k}^{(k)}}
    P_{\vec{n}-\vec{e}_k}(x) Q_{\vec{n} + \vec{e}_k}(y)
\end{equation}
where
\begin{equation} \label{hnk0}
    h_{\vec{n}}^{(k)} = \int_{-\infty}^{\infty} P_{\vec{n}}(x) x^{n_k} w_k(x) dx
\end{equation}
and $\vec{e}_k$ is the $k$th standard basis vector in $\mathbb R^m$.

In the following two sections we study two special cases related to
multiple Hermite polynomials and multiple Laguerre polynomials.
These cases correspond to the random matrix model (\ref{ensemble})
with $V(M) = \frac{1}{2}M^2$ and $V(M) = M$
respectively (in the latter case we restrict to positive
definite matrices). Br\'ezin and Hikami \cite{BH2} and Baik, Ben Arous,
and P\'ech\'e \cite{BBAP} gave double integral representations for
the correlation kernels for these cases. We will derive integral
representations for the multiple Hermite and multiple Laguerre
polynomials, and use that to show that the kernels agree with
the multiple orthogonal polynomial kernel (\ref{CDkernel}).

\section{Multiple Hermite polynomials}

The special case $V(M) = \frac{1}{2} M^2$ was considered in a series
of papers of Br\'ezin and Hikami, \cite{BH1,BH2,BH3,BH4,BH5}. This case corresponds to
$M = H + A$ where $H$ is a random matrix from the GUE ensemble
$ (1/Z_n) e^{-\frac{1}{2} \Tr H^2} dH$
and $A$ is  fixed as before. In \cite{BH2} the following expression
for the kernel was derived
\begin{equation} \label{BHkernel}
    K(x,y) = \frac{1}{(2\pi i)^2}
    \int_{-i\infty}^{i \infty} ds  \oint_{\Gamma} dt \
        e^{ \frac{1}{2} (s-x)^2 - \frac{1}{2} (t-y)^2 } \prod_{k=1}^m
        \left(\frac{s-a_k}{t-a_k}\right)^{n_k} \frac{1}{s-t}
\end{equation}
where $\Gamma$ is a closed contour encircling the points
$a_1, \ldots, a_m$ once in the positive direction, and the
path from $-i\infty$ to $i \infty$ does not intersect $\Gamma$,
see also Johansson \cite{Joh}.

When $V(x) = \frac{1}{2} x^2$, the multiple orthogonal polynomials
are called multiple Hermite polynomials, since they clearly
generalize the usual Hermite polynomials \cite{Apt,ABV,VAC}. We derive integral
representations for the multiple Hermite polynomials of
type I and type II, which resemble the integral representation (\ref{BHkernel})
of the kernel.

\subsection{Multiple Hermite polynomials of type II}

The multiple Hermite polynomial $P_{\vec{n}}$ is the monic
polynomial of degree $|\vec{n}|$ that satisfies (\ref{typeII})
with $w_k(x) = e^{-\frac{1}{2}x^2+a_kx}$, $k=1, \ldots, m$.

\begin{theorem} The multiple Hermite polynomials of type II has the
integral representation
\begin{equation} \label{intrepH2}
    P_{\vec{n}}(x) = \frac{1}{\sqrt{2\pi}i}
    \int_{-i\infty}^{i \infty} e^{\frac{1}{2}(s-x)^2}
        \prod_{k=1}^m (s-a_k)^{n_k} \, ds
\end{equation}
\end{theorem}
\begin{proof}
Let us denote the left-hand side of (\ref{intrepH2}) by $P(x)$.
After performing the change of variables $s = t+x$, we get
\begin{equation} \label{intrepH2b}
    P(x) = \frac{1}{\sqrt{2\pi} i} \int_{-i\infty}^{i \infty} e^{\frac{1}{2} t^2}
    \prod_{k=1}^m (t+x-a_k)^{n_k} \, dt,
\end{equation}
which shows that $P$ is a polynomial of degree $|\vec{n}|$
with leading coefficient
\[ \frac{1}{\sqrt{2\pi} i} \int_{-i\infty}^{i \infty} e^{\frac{1}{2} t^2} dt = 1. \]
So $P$ is a monic polynomial.

Now we use (\ref{intrepH2b}) to compute for $k=1, \ldots, m$,
and $j = 0, 1, \ldots$,
\begin{eqnarray} \nonumber
    \lefteqn{ \int_{-\infty}^{\infty} P(x) x^j e^{-\frac{1}{2} x^2 + a_k x} dx
     = } \\
     && \label{intrepH2c}
     \frac{e^{\frac{1}{2} a_k^2}}{\sqrt{2\pi} i} \int_{-i \infty}^{i\infty} \int_{-\infty}^{\infty}
     e^{\frac{1}{2} (t^2 - (x-a_k)^2)} x^j \prod_{l=1}^m (t + x-a_l)^{n_l} \, dx dt.
\end{eqnarray}
Switching to polar coordinates
$x-a_k = r \cos \theta$, $t = ir \sin \theta$, we find
that the right-hand side of (\ref{intrepH2c}) is equal to
\begin{equation} \label{intrepH2d}
     \frac{e^{\frac{1}{2} a_k^2}}{\sqrt{2\pi}}
     \int_0^{\infty} e^{-\frac{1}{2} r^2} r^{n_k+1}
     \left[\int_0^{2\pi} (a_k + r\cos \theta)^j  e^{i n_k \theta}
        \prod_{l \neq k} (r e^{i \theta} + a_k - a_l)^{n_l} d\theta \right] dr.
\end{equation}
The $\theta$-integral vanishes for $j =0, \ldots, n_k-1$, since the integrand
can be written as a linear combination of $e^{ip\theta}$ with integer $p \geq n_k-j$.
Hence $P$ is the multiple Hermite polynomial of type II and (\ref{intrepH2}) follows.
\end{proof}

\begin{remark}
Evaluating (\ref{intrepH2d}) for $j= n_k$ we find that the $\theta$-integral
is $2\pi \left(\frac{r}{2} \right) \prod_{l \neq k} (a_k-a_l)^{n_l}$ and
\begin{eqnarray} \nonumber
    h_{\vec{n}}^{(k)} & = &
        \int_{-\infty}^{\infty} P_{\vec{n}}(x) x^{n_k} e^{-\frac{1}{2}x^2+a_kx} dx \\
        & = & \sqrt{2\pi} (n_k)! e^{\frac{1}{2} a_k^2} \prod_{l \neq k} (a_k-a_l)^{n_l}.
    \label{hnk}
\end{eqnarray}
\end{remark}

\subsection{Multiple Hermite polynomials of type I}
The multiple Hermite polynomials are polynomials $A_{\vec{n}}^{(k)}$
as in (\ref{degreeAk}) such that the linear form
\begin{equation} \label{QfunctionH}
    Q_{\vec{n}}(x) = \sum_{k=1}^m A_{\vec{n}}^{(k)}(x) e^{-\frac{1}{2} x^2 + a_k x}
\end{equation}
satisfies (\ref{typeI}).

\begin{theorem} The multiple Hermite polynomials of type I have the
integral representation
\begin{equation} \label{intrepH1}
    A_{\vec{n}}^{(k)}(x) e^{-\frac{1}{2} x^2 + a_kx} = \frac{1}{\sqrt{2\pi} 2\pi i}
    \oint_{\Gamma_k} e^{-\frac{1}{2}(t-x)^2}
        \prod_{l=1}^m (t-a_l)^{-n_l} \, dt
\end{equation}
where $\Gamma_k$ is a closed contour encircling $a_k$ once in the positive direction,
but not enclosing any of the other points $a_l$, $l\neq k$.

In addition the function $Q_{\vec{n}}$ of {\rm (\ref{QfunctionH})} has
the integral representation
\begin{equation} \label{intrepH1b}
    Q_{\vec{n}}(x) = \frac{1}{\sqrt{2\pi} 2\pi i}
    \oint_{\Gamma} e^{-\frac{1}{2}(t-x)^2}
        \prod_{l=1}^m (t-a_l)^{-n_l} \, dt
\end{equation}
where $\Gamma$ is a closed contour encircling $a_1, \ldots, a_m$ once in the
positive direction.
\end{theorem}
\begin{proof}
By the residue theorem, we have that the right-hand side of
(\ref{intrepH1}) is equal to
\begin{equation} \label{H1proof1}
    \frac{1}{\sqrt{2\pi}} \frac{1}{(n_k-1)!} \left. \left(\frac{d}{dt}\right)^{n_k-1}
    \left[e^{-\frac{1}{2}(t-x)^2} \prod_{l\neq k} (t-a_l)^{-n_l} \right] \right|_{t=a_k}.
\end{equation}
It is easy to see that (\ref{H1proof1}) has the
form $A_k(x) e^{-\frac{1}{2} x^2 + a_kx}$ where $A_k$ is a polynomial of degree $n_k-1$.
Define the linear form
\begin{equation} \label{H1proof2}
    Q(x) = \sum_{k=1}^m A_k(x) e^{-\frac{1}{2} x^2 + a_kx}
    = \frac{1}{\sqrt{2\pi} 2 \pi i} \oint_{\Gamma} e^{-\frac{1}{2} (t-x)^2}
        \prod_{l=1}^m (t-a_l)^{-n_l} \, dt
\end{equation}
where $\Gamma$ encloses all the points $a_j$, $j=1, \ldots, m$, once
in the positive direction.

Then
\[ \int_{-\infty}^{\infty} x^j Q(x) dx =
    \frac{1}{2\pi i} \oint_{\Gamma}
        \left( \frac{1}{\sqrt{2\pi}} \int_{-\infty}^{\infty}
            x^j e^{-\frac{1}{2}(t-x)^2} dx \right)
            \prod_{l=1}^m (t-a_l)^{-n_l} \, dt.
\]
Since
\[ \frac{1}{\sqrt{2\pi}} \int_{-\infty}^{\infty} x^j e^{-\frac{1}{2} (t-x)^2}  dx =
    \frac{1}{\sqrt{2\pi}} \int_{-\infty}^{\infty} (y+t)^j e^{-\frac{1}{2} y^2} dy
\]
    is a monic polynomial of degree $j$ in the variable $t$,
we find
\begin{equation} \label{H1proof3}
    \int_{-\infty}^{\infty} x^j Q(x) dx =
    \frac{1}{2\pi i} \oint_{\Gamma}
        \pi_j(t) \prod_{l=1}^m (t-a_l)^{-n_l} \, dt, \end{equation}
where $\pi_j$ is a monic polynomial of degree $j$.
Deforming the contour $\Gamma$ to infinity, and using the fact
that the integrand is $s^{j-|\vec{n}|} + O(s^{j-|\vec{n}| - 1})$
as $s \to \infty$, we find that
\[ \int_{-\infty}^{\infty} x^j Q(x) dx =
    0 \qquad \mbox{for } j=0, \ldots, |\vec{n}|-2 \]
and
\[ \int_{-\infty}^{\infty} x^j Q(x) dx =
    1 \qquad \mbox{for } j=|\vec{n}|-1. \]
This shows that $Q = Q_{\vec{n}}$ and $A_k = A_{\vec{n}}^{(k)}$
so that (\ref{intrepH1}) and (\ref{intrepH1b}) follow.
\end{proof}

\subsection{The multiple Hermite kernel}

Let us now show that the Br\'ezin Hikami kernel (\ref{BHkernel}) agrees
with the multiple Hermite kernel (\ref{CDkernel}).
To that end we compute $\frac{\partial K}{\partial x} + \frac{\partial K}{\partial y}$
for the kernel (\ref{BHkernel}) in two ways.

First we have
\begin{eqnarray*}
    \lefteqn{ \frac{\partial K}{\partial x} + \frac{\partial K}{\partial y}
    } \\
    && = \frac{1}{(2\pi i)^2}
    \int_{-i\infty}^{i \infty} ds  \oint_{\Gamma} dt \
        e^{ \frac{1}{2} (s-x)^2 - \frac{1}{2} (t-y)^2 } \prod_{k=1}^m
        \left(\frac{s-a_k}{t-a_k}\right)^{n_k} \frac{-s+x+t-y}{s-t} \\
    &&= (x-y) K(x,y) - \frac{1}{(2\pi i)^2}
    \int_{-i\infty}^{i \infty} ds  \oint_{\Gamma} dt \
        e^{ \frac{1}{2} (s-x)^2 - \frac{1}{2} (t-y)^2 } \prod_{k=1}^m
        \left(\frac{s-a_k}{t-a_k}\right)^{n_k}
\end{eqnarray*}
The last double integral factors into a product of two single integrals,
which by (\ref{intrepH2}) and (\ref{intrepH1b}) leads to
\begin{equation} \label{Kx+Ky1}
    \frac{\partial K}{\partial x} + \frac{\partial K}{\partial y}
    = (x-y) K(x,y) - P_{\vec{n}}(x)Q_{\vec{n}}(y).
\end{equation}

For the second way we evaluate $\frac{\partial K}{\partial x}$ by
noting that $\partial_x e^{\frac{1}{2}(s-x)^2} = - \partial_s e^{\frac{1}{2}(s-x)^2}$,
and integrating by parts the $s$-integral
\begin{eqnarray} \nonumber
    \frac{\partial K}{\partial x}
    & = &
    -\frac{1}{(2\pi i)^2}
    \int_{-i\infty}^{i \infty} ds  \oint_{\Gamma} dt \
        \frac{\partial}{\partial s} e^{ \frac{1}{2} (s-x)^2 -\frac{1}{2} (t-y)^2 }
        \prod_{k=1}^m
        \left(\frac{s-a_k}{t-a_k}\right)^{n_k} \frac{1}{s-t} \\
    &= & \nonumber
    \frac{1}{(2\pi i)^2}
    \int_{-i\infty}^{i \infty} ds  \oint_{\Gamma} dt \
         e^{ \frac{1}{2} (s-x)^2- \frac{1}{2} (t-y)^2 }  \\
         && \qquad \qquad \times \label{Kx}
        \prod_{k=1}^m
        \left(\frac{s-a_k}{t-a_k}\right)^{n_k} \frac{1}{s-t}
        \left\{ \sum_{k=1}^m \frac{n_k}{s-a_k} - \frac{1}{s-t} \right\}.
    \end{eqnarray}
Similarly, we use $\partial_y e^{-\frac{1}{2}(t-y)^2} = - \partial_t e^{-\frac{1}{2}(t-y)^2}$,
and apply integration by parts to the $t$-integral, to obtain
\begin{eqnarray} \nonumber
\frac{\partial K}{\partial y} & = &
    \frac{1}{(2\pi i)^2}
    \int_{-i\infty}^{i \infty} ds  \oint_{\Gamma} dt \
         e^{ \frac{1}{2} (s-x)^2 - \frac{1}{2} (t-y)^2 }  \\
         && \qquad \qquad \times \label{Ky}
        \prod_{k=1}^m
        \left(\frac{s-a_k}{t-a_k}\right)^{n_k} \frac{1}{s-t}
        \left\{ -\sum_{k=1}^m \frac{n_k}{t-a_k} + \frac{1}{s-t} \right\}.
\end{eqnarray}
We add (\ref{Kx}) and (\ref{Ky})
\begin{eqnarray} \nonumber
\frac{\partial K}{\partial x} + \frac{\partial K}{\partial y} & = &
   \frac{1}{(2\pi i)^2}
    \int_{-i\infty}^{i \infty} ds \oint_{\Gamma} dt \
         e^{ \frac{1}{2} (s-x)^2 - \frac{1}{2} (t-y)^2 }  \\
&& \qquad \qquad \times \nonumber
        \prod_{j=1}^m
        \left(\frac{s-a_j}{t-a_j}\right)^{n_j} \frac{1}{s-t}
        \sum_{k=1}^m \left(\frac{n_k}{s-a_k} - \frac{n_k}{t-a_k}\right) \\
& = & - \sum_{k=1}^m n_k \label{Kx+Ky}
    \frac{1}{(2\pi i)^2}
    \int_{-i\infty}^{i \infty} ds \oint_{\Gamma} dt \
         e^{ \frac{1}{2} (s-x)^2 - \frac{1}{2} (t-y)^2 }  \\
         && \qquad \qquad \times   \nonumber
        \prod_{j \neq k} \left(\frac{s-a_j}{t-a_j}\right)^{n_j}
        \frac{(s-a_k)^{n_k-1}}{(t-a_k)^{n_k+1}}.
\end{eqnarray}
For every $k$, the double integral in (\ref{Kx+Ky}) factors into
a product of two single integrals, which by (\ref{intrepH2})
and (\ref{intrepH1b}) are given in terms of multiple Hermite
polynomials. It leads to
\begin{equation} \label{Kx+Ky2}
\frac{\partial K}{\partial x} + \frac{\partial K}{\partial y}
    = - \sum_{k=1}^m n_k P_{\vec{n}-\vec{e}_k}(x) Q_{\vec{n}+\vec{e}_k}(y).
\end{equation}

From (\ref{Kx+Ky1}) and (\ref{Kx+Ky2}) we get
\[ (x-y) K(x,y) = P_{\vec{n}}(x) Q_{\vec{n}}(y) -
    \sum_{k=1}^m n_k P_{\vec{n}-\vec{e}_k}(x) Q_{\vec{n}+\vec{e}_k}(y),
\]
which agrees with (\ref{CDkernel}) since
\begin{equation} \label{nk}
    n_k = \frac{h_{\vec{n}}^{(k)}}{h_{\vec{n}-\vec{e}_k}^{(k)}}
\end{equation}
because of (\ref{hnk}).

\section{Multiple Laguerre polynomials}

Complex Gaussian sample covariance matrices have a distribution
\begin{equation} \label{Wishart}
    \frac{1}{Z_n} e^{- \Tr (\Sigma^{-1} M)} \left(\det M\right)^p dM
\end{equation}
defined on $n \times n$ positive definite Hermitian matrices $M$.
The matrix $M$ arises as $M = X X^H$ where $X$ is an $n \times (n+p)$,
matrix whose independent columns are Gaussian distributed with covariance
matrix $\Sigma$. Here $p$ is a non-negative integer.
The distribution (\ref{Wishart}) is also called a Wishart ensemble.
We assume that $\Sigma^{-1}$ has eigenvalues
$\beta_1, \ldots, \beta_m > 0$ with respective multiplicities $n_1, \ldots, n_m$.

Writing $\Sigma^{-1} = I - A$, we see that (\ref{Wishart}) takes
the form (\ref{ensemble}) with $V(M) = M$, but restricted to positive
definite Hermitian matrices.
It follows that the ensemble (\ref{Wishart}) can be described with
multiple orthogonal polynomials, which in this case are multiple
Laguerre polynomials \cite{Apt,ABV,NS,VAC}. (To be precise, they are
called multiple Laguerre II  in \cite{ABV,VAC} to distinguish them
from another generalization of Laguerre polynomials, called
multiple Laguerre I.)

The eigenvalues of $M$ follow a determinantal
point process on $(0,\infty)$ with kernel $K(x,y)$
given by (\ref{CDkernel})
\begin{equation} \label{CDkernel2}
    (x-y) K(x,y) = P_{\vec{n}}(x) Q_{\vec{n}}(y)
    - \sum_{k=1}^m \frac{h_{\vec{n}}^{(k)}}{h_{\vec{n}-\vec{e}_k}^{(k)}}
        P_{\vec{n}-\vec{e}_k}(x) Q_{\vec{n} + \vec{e}_k}(y)
\end{equation}
where now $P_{\vec{n}}$ is the type II multiple Laguerre polynomial
and $Q_{\vec{n}}(x) = \sum\limits_{k=1}^m A_{\vec{n}}^{(k)}(x) x^p e^{-\beta_k}(x)$
is the linear form involving the type I multiple Laguerre polynomials
$A_{\vec{n}}^{(k)}$.

Baik, Ben Arous and P\'ech\'e \cite{BBAP} gave a double integral representation
for the correlation kernel
\begin{equation} \label{BBAPkernel}
    K(x,y) = \frac{1}{(2\pi i)^2} \oint_{\Sigma} ds \oint_{\Gamma} dt
    e^{xs-yt} \left(\frac{t}{s}\right)^{|\vec{n}|+p}
        \prod_{k=1}^m \left( \frac{s-\beta_k}{t-\beta_k}
        \right)^{n_k} \frac{1}{s-t}
\end{equation}
where $\Sigma$ and $\Gamma$ are disjoint closed contours both oriented counterclockwise
such that $\Sigma$ encloses $0$ and lies in $\{s \in \mathbb C \mid \Re s < \min_k \beta_k \}$
and $\Gamma$ encloses the points $\beta_1, \ldots, \beta_m$ and lies in the right half-plane.

In view of (\ref{BBAPkernel}) and our experience with multiple Hermite polynomials
we expect integral representations for the multiple Laguerre polynomials as well.
We will see that this is indeed the case, and we use this to study the connection
between the kernels (\ref{CDkernel2}) and (\ref{BBAPkernel}). It will turn
out that the two kernels are equal up to a multiplicative factor $x^p y^{-p}$.
However, this difference does not affect the correlation functions (\ref{correlations}).

\subsection{Multiple Laguerre polynomials of type II}
The multiple Laguerre polynomial of type II is a monic polynomial $P_{\vec{n}}$
of degree $|\vec{n}|$ such that
\begin{equation} \label{LaguerreII}
    \int_0^{\infty} P_{\vec{n}}(x) x^{j+p} e^{-\beta_k x} dx = 0,
    \qquad k=1, \ldots, m, \quad j=0, \ldots, n_k-1.
\end{equation}

\begin{theorem}
    The multiple Laguerre polynomial of type II has the integral representation
    \begin{equation} \label{intrepL2}
        P_{\vec{n}}(x) =
        \frac{(|\vec{n}|+p)! x^{-p}}{2 \pi i \prod_{k=1}^m (-\beta_k)^{n_k}}
         \oint_{\Sigma} e^{xs} s^{-|\vec{n}|-p-1}
            \prod_{k=1}^m (s- \beta_k)^{n_k} ds
    \end{equation}
    where $\Sigma$ is a closed contour around $0$ oriented counterclockwise,
    which does not enclose any of the $\beta_k$'s.
\end{theorem}
\begin{proof}
    Denote the right-hand side of (\ref{intrepL2}) by $P(x)$
    and write $C = \frac{(|\vec{n}|+p)!}{2 \pi i \prod_{k=1}^m (-\beta_k)^{n_k}}$.
    Then
    \begin{equation} \label{intrepL2b}
    P(x) = C x^{-p} \sum_{j=0}^{\infty} \frac{x^j}{j!}
    \oint_{\Sigma} s^{j-|\vec{n}|-p-1} \prod_{k=1}^m (s-\beta_k)^{n_k} ds.
    \end{equation}
    All terms in (\ref{intrepL2b}) with $j \geq |\vec{n}| + 1$ vanish by Cauchy's theorem,
    as well as all terms with $j \leq p-1$, as we can see by deforming the contour
    $\Sigma$ to infinity. It follows that $P$ is a polynomial of degree $|\vec{n}|$,
    whose leading coefficient is
    \[ C \frac{1}{(|\vec{n}|+p)!}
    \oint_{\Sigma} s^{-1} \prod_{k=1}^m (s-\beta_k)^{n_k} ds
    = C \frac{2\pi i}{(|\vec{n}|+p)!} \prod_{k=1}^m (- \beta_k)^{n_k} = 1. \]
    Hence $P$ is a monic polynomial of degree $|\vec{n}|$.

    We now verify the orthogonality conditions (\ref{LaguerreII}).
    Take $k = 1, \ldots, m$ and let $j = 0, \ldots, n_k-1$. Then
    \[ \int_0^{\infty} P(x) x^{j+p} e^{-\beta_k x} dx =
        C \oint_{\Sigma}
        \left( \int_0^{\infty} x^j e^{(s-\beta_k)x} dx \right) s^{-|\vec{n}|-p-1}
            \prod_{l=1}^m (s- \beta_l)^{n_l} ds \]
    where we have assumed that $\Sigma$ is so that $\Re s <  \beta_k$
    for every $s \in \Sigma$.
    The inner integral is  $j! (\beta_k-s)^{-j-1}$, so that
    \[ \int_0^{\infty} P(x) x^{j+p} e^{-\beta_kx} dx =
        C j! \oint_{\Sigma} s^{-|\vec{n}|-p-1} (s-\beta_k)^{n_k-j-1} \prod_{l \neq k} (s-\beta_l)^{n_l} ds \]
        The integrand in the last integral has no pole at $\beta_k$, since
         $j \leq n_k-1$. At infinity the integrand behaves like $s^{-j-p-2}$. So by deforming the
         contour $\Sigma$ to infinity, we conclude that the integral is zero.
         The conclusion is that $P$ is the multiple Laguerre polynomial
        of type II and (\ref{intrepL2}) follows.
\end{proof}

\subsection{Multiple Laguerre polynomials of type I}
The multiple Laguerre polynomials of type I $A_{\vec{n}}^{(k)}$, $k=1, \ldots, m$, have
degrees
\begin{equation} \label{LaguerreIa}
    \deg A_{\vec{n}}^{(k)} \leq n_k-1, \qquad k =1,\ldots, m
\end{equation}
and are such that
\[ Q_{\vec{n}}(x) = \sum_{k=1}^m A_{\vec{n}}^{(k)} x^p e^{-\beta_kx} \]
satisfies
\begin{equation} \label{LaguerreIb}
    \int_0^{\infty} x^j Q_{\vec{n}}(x) dx =
    \left\{ \begin{array}{ll} 0, & \quad j=0, \ldots, |\vec{n}|-2, \\
    1, & \quad j=|\vec{n}|-1.
    \end{array} \right.
\end{equation}
There is also an integral representation for the  multiple Laguerre polynomials
of type I.

\begin{theorem} The multiple Laguerre polynomials of type I have the
integral representation
\begin{equation} \label{intrepL1}
    A_{\vec{n}}^{(k)}(x) x^p e^{-\beta_k x} = -
        \frac{\prod_{l=1}^m (-\beta_l)^{n_l} x^p}{2 \pi i (|\vec{n}|+p-1)! }
         \oint_{\Gamma_k} e^{-xt} t^{|\vec{n}|+p-1}
            \prod_{l=1}^m (t- \beta_l)^{-n_l} dt
\end{equation}
where $\Gamma_k$ is a closed contour around $\beta_k$, which does not enclose $0$
nor any of the other points $\beta_l$, $l\neq k$.

In addition the function $Q_{\vec{n}}$ has the integral representation
\begin{equation} \label{intrepL1b}
    Q_{\vec{n}}(x) = -
        \frac{\prod_{l=1}^m (-\beta_l)^{n_l} x^p}{2 \pi i (|\vec{n}|+p-1)! }
         \oint_{\Gamma} e^{-xt} t^{|\vec{n}|+p-1}
            \prod_{l=1}^m (t- \beta_l)^{-n_l} dt
\end{equation}
where $\Gamma$ is a closed contour around $\beta_1, \ldots, \beta_m$, but which
does not enclose $0$.
\end{theorem}
\begin{proof}
Only the pole at $\beta_k$ contributes to the integral in (\ref{intrepL1}).
By the residue theorem, we have that the right-hand side of (\ref{intrepL1})
\[ const  \ x^p \left(\frac{d}{dt}\right)^{n_k-1}
    \left. \left[ e^{-xt} t^{|\vec{n}|+p-1} \prod_{l\neq k} (t-\beta_l)^{-n_l} \right]
        \right|_{t=\beta_k}
\]
which is easily seen to of the form $A_k(x) x^p e^{-\beta_k x}$
where $A_k$ is a polynomial of degree $n_k-1$.

Let \[ Q(x) = \sum_{k=1}^m A_k(x) x^p e^{-\beta_k x} \]
which is then equal to the right-hand side of (\ref{intrepL1b}).
We may assume that $\Gamma$ is entirely in the right half-plane.
Then
\[ \int_0^{\infty} x^j Q(x) dx = -
        \frac{\prod_{l=1}^m (-\beta_l)^{n_l}}{2 \pi i (|\vec{n}|+p-1)! }
         \oint_{\Gamma} dt \int_0^{\infty} dx x^{j+p} e^{-xt} \ t^{|\vec{n}|+p-1}
            \prod_{l=1}^m (t- \beta_l)^{-n_l}. \]
The $x$-integral is $(j+p)! t^{-j-p-1}$, so that
\begin{equation} \label{intrepL1c}
    \int_0^{\infty} x^j Q(x) dx = -
        \frac{\prod_{l=1}^m (-\beta_l)^{n_l} (j+p)!}{2 \pi i (|\vec{n}|+p-1)! }
         \oint_{\Gamma} t^{|\vec{n}|-j-2}
            \prod_{l=1}^m (t- \beta_l)^{-n_l} dt.
\end{equation}
Assuming $j \leq |\vec{n}|-2$, we can deform $\Gamma$ to infinity
without picking up a residue contribution at $t=0$.
The integrand behaves like $t^{-j-2}$ at infinity, and so (\ref{intrepL1c})
vanishes for $j \leq |\vec{n}|-2$. For $j= |\vec{n}|-1$, we pick up
a residue contribution at $t=0$, and the result is that
\begin{equation} \label{intrepL1d}
\int_0^{\infty} x^{|\vec{n}|-1} Q(x) dx = 1.
\end{equation}
Thus $Q$ satisfies the type I multiple orthogonality conditions (\ref{LaguerreIb})
and the theorem follows.
\end{proof}

\subsection{The multiple Laguerre kernel}
We finally compare the representations (\ref{CDkernel2}) and (\ref{BBAPkernel})
of the kernel.

We start from the double integral (\ref{BBAPkernel}) and evaluate $xK(x,y)$
using an integration by parts on the $s$-integral. The result is
\begin{eqnarray} \nonumber
 xK(x,y) & = & \frac{1}{(2\pi i)^2} \oint_{\Sigma} ds \oint_{\Gamma} dt
    e^{xs-yt} \left(\frac{t}{s}\right)^{|\vec{n}|+p}
        \prod_{l=1}^m \left( \frac{s-\beta_l}{t-\beta_l}
        \right)^{n_l}  \\
        & & \qquad \times \label{xKxy} \ \frac{1}{s-t}
        \left\{ \frac{|\vec{n}|+p}{s} - \sum_{k=1}^m \frac{n_k}{s-\beta_k} + \frac{1}{s-t} \right\}.
\end{eqnarray}
Similarly, after integration by parts on the $t$-integral,
\begin{eqnarray} \nonumber
 yK(x,y) & = & \frac{1}{(2\pi i)^2} \oint_{\Sigma} ds \oint_{\Gamma} dt
    e^{xs-yt} \left(\frac{t}{s}\right)^{|\vec{n}|+p}
        \prod_{l=1}^m \left( \frac{s-\beta_l}{t-\beta_l}
        \right)^{n_l}  \\
        & & \qquad \times \label{yKxy} \ \frac{1}{s-t}
        \left\{ \frac{|\vec{n}|+p}{t} - \sum_{k=1}^m \frac{n_k}{t-\beta_k}  + \frac{1}{s-t} \right\}.
\end{eqnarray}
Hence
\begin{eqnarray*}
\lefteqn{(x-y) K(x,y) } \\
&& = \frac{1}{(2\pi i)^2} \oint_{\Sigma} ds \oint_{\Gamma} dt
    e^{xs-yt} \left(\frac{t}{s}\right)^{|\vec{n}|+p}
        \prod_{l=1}^m \left( \frac{s-\beta_l}{t-\beta_l}
        \right)^{n_l} \\
        & & \qquad \times \
        \left\{ - \frac{|\vec{n}|+p}{st} + \sum_{k=1}^m \frac{n_k}{(s-\beta_k)(t-\beta_k)} \right\} \\
        & & = -\frac{|\vec{n}|+p}{(2\pi i)^2} \oint_{\Sigma} ds \oint_{\Gamma} dt
    e^{xs-yt} \frac{t^{|\vec{n}|+p-1}}{s^{|\vec{n}|+p+1}}
    \prod_{l=1}^m \left( \frac{s-\beta_l}{t-\beta_l}  \right)^{n_l} \\
    & & \quad + \sum_{k=1}^m \frac{n_k}{(2\pi i)^2} \oint_{\Sigma} ds \oint_{\Gamma} dt
    e^{xs-yt} \left(\frac{t}{s}\right)^{|\vec{n}|+p}
        \prod_{l\neq k} \left( \frac{s-\beta_l}{t-\beta_l}
        \right)^{n_l} \frac{(s-\beta_k)^{n_k-1}}{(t-\beta_k)^{n_k+1}}
\end{eqnarray*}
Now we have $m+1$ double integrals and they all factor into products
of two single integrals of the forms (\ref{intrepL2}) and (\ref{intrepL1b}).
The result is that
\begin{equation} \label{compareK}
    (x-y) K(x,y) =
    x^p y^{-p} \left( P_{\vec{n}}(x) Q_{\vec{n}}(y) -
        \sum_{k=1}^m n_k \frac{|\vec{n}|+p}{\beta_k^2}  P_{\vec{n}-\vec{e}_k}(x) Q_{\vec{n}+\vec{e}_k}(y)
        \right).
\end{equation}
It can be shown that
\[ n_k \frac{|\vec{n}|+p}{\beta_k^2} = \frac{h_{\vec{n}}^{(k)}}{h_{\vec{n}-\vec{e}_k}^{(k)}} \]
so that (\ref{compareK}) agrees with (\ref{CDkernel2}) up to the factor $x^p y^{-p}$.
However, this factor is not essential since it does not change the correlation functions
(\ref{correlations}). Hence (\ref{CDkernel2}) and (\ref{BBAPkernel}) are essentially the same.


\begin{thebibliography}{99}
\bibitem{Apt}
    A. I. Aptekarev, Multiple orthogonal polynomials,
    J. Comput. Appl. Math. 99 (1998), 423--447.
\bibitem{ABV}
    A. I. Aptekarev, A. Branquinho, and W. Van Assche,
    Multiple orthogonal polynomials for classical weights,
    Trans. Amer. Math. Soc. 355 (2003), 3887--3914.
\bibitem{BBAP}
    J. Baik, G. Ben Arous, and S. P\'ech\'{e},
    Phase transition of the largest eigenvalue for non-null
    complex sample covariance matrices, preprint math.PR/0403022.
\bibitem{BK1} P. M. Bleher and A. B. J. Kuijlaars,
    Random matrices with external source and multiple orthogonal polynomials,
    Internat. Math. Research Notices 2004 (2004), 109--129.
\bibitem{BK2} P. M. Bleher and A. B. J. Kuijlaars,
    Large $n$ limit of Gaussian random matrices with external source, part I,
    preprint math-ph/0402042.
\bibitem{Bor} A. Borodin, Biorthogonal ensembles,
    Nuclear Phys. B 536 (1999),  704--732.
\bibitem{BH1} E. Br\'ezin and S. Hikami,
    Correlations of nearby levels induced by a random potential,
    Nucl. Phys. B 479 (1996), 697--706.
\bibitem{BH2} E. Br\'ezin and S. Hikami,
    Spectral form factor in a random matrix theory,
    Phys. Rev. E 55 (1997), 4067--4083.
\bibitem{BH3} E. Br\'ezin and S. Hikami,
    Extension of level-spacing universality,
    Phys. Rev. E 56 (1997), 264--269.
\bibitem{BH4} E. Br\'ezin and S. Hikami,
    Universal singularity at the closure of a gap in a random matrix theory,
    Phys. Rev. E 57 (1998), 4140--4149.
\bibitem{BH5} E. Br\'ezin and S. Hikami,
    Level spacing of random matrices in an external source,
    Phys. Rev. E 58 (1998), 7176--7185.
\bibitem{DK} E. Daems and A. B. J. Kuijlaars,
    A Christoffel-Darboux formula for multiple orthogonal polynomials,
    preprint math.CA/0402031.
\bibitem{HC} Harish-Chandra, Differential operators on a
    semisimple Lie algebra, Amer. J. Math. 79 (1957), 87--120.
\bibitem{IZ} C. Itzykson and J. B. Zuber, The planar
    approximation II, J. Math. Phys. 21 (1980), 411--421.
\bibitem{Joh} K. Johansson,
    Universality of the local spacing distribution in certain ensembles
    of Hermitian Wigner matrices,
    Comm. Math. Phys. 215 (2001), no. 3, 683--705.
\bibitem{NS} E. Nikishin and V. Sorokin,
    Rational Approximation and Orthogonality,
    Translations of Mathematical  Monographs 92, Amer. Math. Soc. Providence R.I., 1991.
\bibitem{VAC} W. Van Assche and E. Coussement,
    Some classical multiple orthogonal polynomials,
    J. Comput. Appl. Math. 127 (2001), 317--347.
\bibitem{ZJ1} P. Zinn-Justin,
    Random Hermitian matrices in an external field,
    Nuclear Phys. B 497 (1997), 725--732.
\bibitem{ZJ2} P. Zinn-Justin,
    Universality of correlation functions of Hermitian random matrices
    in an external field,
    Comm. Math. Phys. 194 (1998), 631--650.
\end{thebibliography}
\end{document}